\documentclass[12pt]{article}
\usepackage{amsmath}
\usepackage{amsfonts}
\usepackage{amssymb, amsthm}

\newcommand{\im}{\mathrm{im}\,}

\newtheorem*{theorem*}{Theorem}
\newtheorem{lemma}{Lemma}
\newtheorem*{claim*}{Claim}
\newtheorem*{def*}{Definition}
\newtheorem*{conj*}{Conjecture}
\newtheorem*{remark*}{Remark}

\begin{document}
\title{New equation on the low dimensional Calabi--Yau metrics\footnote{This work was supported in part by Russian Foundation for Basic
Research (grant 09-01-00598-a) and the Council of the Russian
Federation Presidential Grants (projects NSh-7256.2010.1 and
MK-842.2011.1).}}
\author{Dmitry  Egorov}
\date{}
\maketitle \sloppy

\begin{abstract}
In this paper we introduce a new equation on the compact K\"{a}hler
manifolds. Solution of this equation corresponds to the Calabi--Yau
metric. New equation differs from the Monge--Amp\`{e}re equation
considered by Calabi and Yau.
\end{abstract}

\section{Introduction}

Let $(M,\omega)$ be the compact K\"{a}hler  $n$-manifold. Calabi
made the following conjecture proven by Yau: If $c_1(M)=0$, then
there exists a Riemannian metric with the holonomy group contained
in $SU(n)$ \cite{Calabi, Yau}. Explicitly, Calabi conjecture is
stated as follows: There exists unique real function $\varphi$ such
that:
\begin{equation}
\label{monge_ampere} (\omega+i\partial\bar{\partial}\varphi)^n =
e^F\omega^n, \quad \int_M{\varphi}\,\omega^n = 0,
\end{equation}
where $F$ is a smooth real function on $M$ such that
$\int_M{e^F}\omega^n = \int_M{\omega^n}$. It is assumed that
$\omega+i\partial\bar{\partial}\varphi >0$.

Existence of the solution for \eqref{monge_ampere} implies existence
of the Ricci-flat metric on $M$.

Note that the Monge--Amp\`{e}re equation describes the deformation
of the K\"{a}hler form. We study the deformation of the holomorphic
volume form on the K\"{a}hler manifold $M$ in complex dimensions $2$
and $3$.

\begin{theorem*}
Let $M$ be a compact K\"{a}hler $n$-manifold. Suppose $c_1(M)=0$.
Denote holomorphic volume form by $\Omega$. For $n=2, 3$ there
exists a solution of the following equation:
\begin{equation}\label{new_equation}
(\Omega+dd^s\psi)\wedge(\bar{\Omega} +dd^s\bar{\psi}) =
e^F\Omega\wedge\bar{\Omega},
\end{equation}
where $\psi$ is a complex $n$-form such that $\tilde{\Omega} =
\Omega+dd^s\psi$ is a stable primitive form; $F$ is a smooth real
function such that $\int_M{e^F\Omega\wedge\bar{\Omega}} = \int_M
\Omega\wedge\bar{\Omega}$. Here $d^s$ is the symplectic differential
operator defined in \S\ref{section_symplectic}.
\end{theorem*}
Here $\tilde{\Omega}$ is said to be  stable if there exists some
complex structure $J$ on $M$ such that $\tilde{\Omega}$ is a
holomorphic volume form with respect to $J$.

The main theorem  implies the existence of the Calabi--Yau metric on
$M$. Since $\tilde{\Omega} = \Omega+dd^s\psi$ is the holomorphic
volume form,  $\tilde{\Omega}\wedge\bar{\tilde{\Omega}}$ is
proportional to the volume form.Then there exists a real function
$F$ such that
$$
e^F{\Omega}\wedge\bar{{\Omega}} = \omega^n,
$$
where $F$ is a smooth real function such that
$\int_M{e^F\Omega\wedge\bar{\Omega}} =
\int_M{\Omega\wedge\bar{\Omega}}$. By the main theorem, there exists
$\tilde{\Omega}$ such that
\begin{equation}
\label{cy_monge_ampere}\tilde{\Omega}\wedge\bar{\tilde{\Omega}} =
{\omega}^n.
\end{equation}
Recall that $n$-form is called primitive if
\begin{equation}
\label{cy_compat}\Omega\wedge{\omega} = 0.
\end{equation}

Stable forms $\tilde{\Omega}$ and $\omega$ are closed and satisfy
\eqref{cy_monge_ampere},\eqref{cy_compat}. Therefore they set up an
integrable $SU(n)$-structure. Hence the holonomy group of $M$ is
contained in $SU(n)$.

In \S \ref{section_symplectic} we recall necessary facts from the
symplectic Hodge theory. In \S \ref{section_proof_1} we prove the
main theorem. The proof follows from Yau's theorem. In \S
\ref{conclusion} we discuss possible generalizations of the main
theorem.

The author is grateful to Iskander Taimanov for drawing attention to
this area of mathematics and constant support and to Andrey Mironov
for useful conservations.

\section{Symplectic Hodge theory}
\label{section_symplectic} Let $(M,\omega)$ be a symplectic manifold
of real dimension $2n$. By $*_s$ denote a symplectic Hodge star. The
action of $*_s$ on the $k$-forms is uniquely determined by the
following formula:
$$
\alpha\wedge *_s\beta =
(\omega^{-1})^k(\alpha,\beta)\frac{\omega^n}{n!}.
$$
By definition, put
\begin{equation}\label{delta} d^s =
(-1)^{k+1}\ast_s d\,\ast_s.
\end{equation}
Note that $d^s$ reduces degree of the form by one and  $dd^s =
-d^sd$.
K\"{a}hler manifolds satisfy the following condition ($dd^s$-lemma):
$$
\im d \cap  \ker d^s = \ker d \cap \im d^s = \im dd^s.
$$

\begin{lemma}
\label{lemma1} Let $M$ be a K\"{a}hler $n$-manifold. Suppose
$\Omega$ and $\tilde{\Omega}$ are primitive cohomologous $n$-forms.
Then $\tilde{\Omega} = \Omega + dd^s\psi$, where $\psi$ is some
$n$-form.
\end{lemma}
\begin{proof}
The following formula for primitive $n$-forms \cite{Weil} (see also
\cite{Yau_Tseng}):
\begin{equation}
\label{star_formula} \ast_s \eta = (-1)^{\frac{n(n+1)}{2}}\eta
\end{equation}
implies that $d$-closed forms are $d^s$-closed as well. By the
$dd^s$-lemma, $\tilde{\Omega}-\Omega$ takes the form $dd^s\psi$.
\end{proof}

\section{Proof of the main theorem}
\label{section_proof_1} Recall the definition of a stable form by
Hitchin \cite{Hitchin3f,Hitchin_stable}. Let $V$ be a real
$m$-dimensional space. A real form $\rho\in \Lambda^pV^*$ is called
stable if the orbit of $\rho$ under the natural action of $GL(V)$ is
open.

If $p=2$ and $m=2n$, then we obtain the non-degeneracy condition of
the symplectic form. Recall that any compact K\"{a}hler $2$-manifold
with $c_1(M)=0$ is hyper-K\"{a}hler. Then by dimensional reasons
real and imaginary parts of the holomorphic volume form are
symplectic and hence are stable.

If $p=3$ and $m=6$, then by definition stable forms are real and
imaginary parts of the holomorphic volume form.

We need the Moser theorem for oriented manifolds \cite{Moser}.
\begin{theorem*}
Let $M$ be an oriented manifold. Suppose there exists a family of
cohomologous symplectic forms $\{\omega_t\}, {t\in [0,1]}$. Then
there exists an orientation preserving diffeomorphism $\phi_t$ such
that
$$
\phi_t^*\omega_t = \omega_0,
$$
and  $\phi_0$ is the identity map.
\end{theorem*}
\begin{proof}
By Yau's theorem, there exists $\tilde{\omega}\sim\omega$ compatible
with $\Omega$.

In his proof Yau uses the continuity method. Therefore we get a
family of cohomologous symplectic forms such that
$$
\omega_0 = \omega, \quad \omega_1 = \tilde{\omega}.
$$
Denote $\phi_1$ by $\phi$. We claim that $\phi^*\Omega$ is the
required $\tilde{\Omega}$. Indeed, act by $\phi^*$ on the compatible
pair: $\Omega$ and $\tilde{\omega} = (\phi^{-1})^*\omega$.
Diffeomorphism preserves compatibility conditions
\eqref{cy_monge_ampere},\eqref{cy_compat} and stability. Since
$\phi$ is isotopic to the identity map, $\tilde{\Omega} =
\phi^*\Omega \in [\Omega]$.

Stable forms $\tilde{\Omega}$ and $\omega$ are closed and compatible
in the sense of \eqref{cy_monge_ampere},\eqref{cy_compat}. Hence
they define a Riemannian metric with holonomy group contained in
$SU(n)$.

Let's prove that $\tilde{\Omega}$ is a solution of
\eqref{new_equation}. Since $\tilde{\Omega}$ is the holomorphic
volume form with respect to some complex structure,
$\tilde{\Omega}\wedge\bar{\tilde{\Omega}}$ is proportional to the
volume form. Therefore there exists a smooth real function $F$ such
that
\begin{equation}
\label{fff}  e^F\tilde{\Omega}\wedge\bar{\tilde{\Omega}} = \omega^n,
\end{equation}
$$
\int_M{e^F\tilde{\Omega}\wedge\bar{\tilde{\Omega}}} =
\int_M{\tilde{\Omega}\wedge\bar{\tilde{\Omega}}}.
$$
Substituting $\eqref{fff}$ in $\eqref{cy_monge_ampere}$, we obtain
the required equation \eqref{new_equation}:
$$
\tilde{\Omega}\wedge\bar{\tilde{\Omega}} = \omega^n =
e^F\Omega\wedge\bar{\Omega}.
$$
\end{proof}

The key moment of the proof is the stability of the holomorphic
volume form. Therefore this proof can not be generalized to higher
dimensions.

\section{Conclusion}
\label{conclusion} Here we discuss possible extensions of the work.

Let $(M,\omega)$ be a symplectic non-K\"{a}hler manifold of
dimension $2n$. Direct generalization of the Monge--Amp\`{e}re
equation is as follows:
\begin{equation}
\label{djd_eqn} (\omega + dJd\varphi)^n = e^F\omega^n,
\end{equation}
where $J$ is an almost  complex structure. If $J$ is integrable,
then this equation turns into the Monge--Amp\`{e}re equation.

Equation \eqref{djd_eqn} doesn't give anything new. It is proved in
\cite{Delanoe,Wang_Zhu} that if for any compatible RHS there exists
a solution, then $J$ is integrable.

It is interesting to find an analogue of this result for new
equation. In generalized complex geometry a non-degenerate $2$-form
is called an almost symplectic structure. Integrability stands for
closedness.

Assume that for any compatible RHS there exists a solution of the
equation:
$$
\tilde{\Omega}\wedge\bar{\tilde{\Omega}} =
e^F\Omega\wedge\bar{\Omega}, \quad \tilde{\Omega} = \Omega +
dd^s\psi.
$$
Does it imply integrability, i.e., closedness, of the almost
symplectic structure?

Let's introduce the second possible extension of the current work.
The Monge--Amp\`{e}re equation is generalized to non-K\"{a}hler
symplectic $4$-manifolds by Donaldson--Yau--Weinkove--Tosatti. It is
suggested in \cite{Donaldson} to consider the following equation:
\begin{equation}
\label{cy_eqn} \tilde{\omega}^n = e^F\omega^n,
\end{equation}
where $\tilde{\omega}\sim\omega$ are symplectic forms. This equation
is called the Calabi--Yau equation.

In \cite{W, TWY, TW} authors assume that the curvature of the Chern
connection satisfies some restrictions. Then they prove the
existence of the unique solution for \eqref{cy_eqn}.

Consider a generalization of the new equation to non-K\"{a}hler
manifolds. In our opinion the non-K\"{a}hler Hermitian manifolds
with trivial canonical bundle (see for example \cite{Fine_Panov})
are natural candidates for generalization of the existence theorem.

Finally, in the K\"{a}hler case it is unclear: Does there exist
generalization of the existence theorem for $n>3$ ?

\vskip5mm {\sc Ammosov Northeastern Federal University, Yakutsk,
677000,
Russia}\\
\textit{e-mail:}{\tt egorov.dima@gmail.com}

\end{document}